\documentclass[12pt]{amsart}
\usepackage{amssymb,times}

\makeatletter
\def\@strippedMR{}
\def\@scanforMR#1#2#3\endscan{
  \ifx#1M\ifx#2R\def\@strippedMR{#3}
  \else\def\@strippedMR{#1#2#3}
  \fi\fi}
\renewcommand\MR[1]{\relax\ifhmode\unskip\spacefactor3000 \space\fi
  \@scanforMR#1\endscan
  MR\MRhref{\@strippedMR}{\@strippedMR}}
\makeatother

\parskip=\medskipamount

\newtheorem*{Thm*}{Theorem}
\newtheorem{Thm}{Theorem}

\newtheorem{Prop}[Thm]{Proposition}

\theoremstyle{definition}

\newtheorem{Remark}{Remark}

\newtheorem{Question}{Question}

\newcommand{\set}[1]{\left\{#1\right\}}

\renewcommand{\phi}{\varphi}

\allowdisplaybreaks[1]

\title{A characterization of ultraspherical, Hermite, and Chebyshev polynomials of the first kind}
\author{Michael Anshelevich}
\thanks{This work was supported in part by the NSF grants DMS-0900935 and DMS-1160849.}
\address{Department of Mathematics, Texas A\&M University, College Station, TX 77843-3368}
\email{manshel@math.tamu.edu}
\date{\today}

\begin{document}

\begin{abstract}
We show that the only orthogonal polynomials with a generating function of the form $F \left(x z - \alpha z^2 \right)$ are the ultraspherical, Hermite, and Chebyshev polynomials of the first kind. For special $F$ for which this is the case, we then finish the classification of orthogonal polynomials with more general generating functions $F(x w(z) - R(z))$.
\end{abstract}

\maketitle

\section{Introduction}

In this short note, we answer the following two questions.

\begin{Question}
\label{Question:Main}
What are all the orthogonal polynomials with generating functions of the form
\[
F \left(x z - \alpha z^2 \right)
\]
for some number $\alpha$ and function (or, more precisely, formal power series) $F$?
\end{Question}

\begin{Question}
For $F$ appearing in the answer to Question~\ref{Question:Main}, what are all the orthogonal polynomials with more general generating functions
\[
F(x w(z) - R(z)),
\]
for formal power series $w, A$ with $w(0) = 1$, $R(0) = 0$?
\end{Question}

The rest of the introduction explains the motivation and context for these questions.

The question of characterizing various classes of orthogonal polynomials has a long and distinguished history, see \cite{Al-Salam-Characterization} for an excellent survey up to 1990. The study of general polynomial families goes back to Paul Appell in 1880 \cite{Appell}, who looked at polynomials with generating functions of the form
\[
\sum_{n=0}^\infty \frac{1}{n!} P_n(x) z^n = A(z) \exp(x z)
\]
for some function $A(z)$. These are now called Appell polynomials. Later, they were generalized to Sheffer families with generating functions
\begin{equation}
\label{Meixner}
\sum_{n=0}^\infty \frac{1}{n!} P_n(x) z^n = A(z) \exp(x U(z))
\end{equation}
for some functions $A(z), U(z)$. The prototypical ``orthogonal polynomials characterization result'' is Meixner's 1934 description of all \emph{orthogonal} polynomials with the Sheffer-type generating functions \cite{Meixner}. On the other hand, among Appell polynomials, only Hermite polynomials are orthogonal.

Besides nice generating functions, the Meixner class has many other characterizations and applications, see \cite{Diaconis-Gibbs-sampling} for an excellent survey. Perhaps for this reason, many generalizations of this class have been attempted. The most popular of these are probably the $q$-deformed families. One approach (there are several) extends the Sheffer class by looking at the generating functions of the form
\[
A(z) \prod_{k=0}^\infty \frac{1}{1 - (1-q) U \left(q^k z \right) z}
\]
(after appropriate normalization, one gets the Sheffer form for $q \rightarrow 1$). In this case the analog of the Meixner class are the Al-Salam and Chihara polynomials \cite{Al-Salam-Pollaczek}. For the study of two different types of $q$-Appell polynomials, see \cite{Al-Salam-q-Appell,Al-Salam-Hermite}.

A different generalization of the Sheffer class are generating functions of the general Boas-Buck \cite{Boas-Buck-Book} type:
\[
\sum_{n=0}^\infty c_n P_n(x) z^n = A(z) F(x U(z))
\]
for $F(z) = \sum_{n=0}^\infty c_n z^n$ with $c_0 = 1$. The usual case corresponds to $F(z) = e^z$. In the Boas-Buck setting, the problem of describing all orthogonal polynomials is wide open. The Appell-type class (with $U(z) = z$) in this case consists of the Brenke polynomials, and at least in that case all the orthogonal polynomials are known \cite{Chihara-Brenke-1,Asai-Kubo-Kuo-Brenke}.

Now note that in the Sheffer/Meixner case in equation~\eqref{Meixner}, corresponding to $F(z) = e^z$, the generating function has an alternative form
\begin{equation*}
\sum_{n=0}^\infty \frac{1}{n!} P_n(x) z^n = A(z) \exp(x U(z)) = \exp(x U(z) + \log A(z)).
\end{equation*}
So another interesting class to look at are (all or just orthogonal) polynomials with generating functions
\begin{equation}
\label{Exp}
\sum_{n=0}^\infty c_n P_n(x) z^n = F(x w(z) - R(z)),
\end{equation}
which again gives the Sheffer/Meixner families for $F(z) = e^z$.

The case $F(z) = \frac{1}{1-z}$ appears in Free Probability \cite{Nica-Speicher-book}, see Section~3 of \cite{AnsMeixner} for the author's description of the ``free Meixner class'', which is in a precise bijection with the Meixner class (except for the binomial case \cite{Boz-Bryc}). In fact, this family was already described in \cite{Al-Salam-Verma}. Here again, one can write the generating function in two ways:
\[
A(z) \frac{1}{1 - x U(z)} = \frac{1}{1 - \left( x \frac{U(z)}{A(z)} - \frac{1 - A(z)}{A(z)} \right)}.
\]
More generally, Boas and Buck proved the following result.
\begin{Thm*}
\cite{Boas-Buck-Article}
The only functions $F$ with $F(0) = 1$ such that
\begin{equation}
\label{Boas-Buck-two-classes}
A(z) F(x U(z)) = F(x w(z) - R(z))
\end{equation}
are $F(z) = e^z$ and $F(z) = \frac{1}{(1 - z)^\lambda}$ for some $\lambda$.
\end{Thm*}

So as an alternative to the Boas-Buck formulation, we are interested in orthogonal polynomials with generating functions of the form $F(x w(z) - R(z))$, or at least in the Appell-type subclass $F(x z - R(z))$. A priori, for general $R$ even this seems to be a hard question. However, the orthogonal Appell polynomials are only the Hermite polynomials, with the exponential generating function
\[
\exp \left(x z - z^2/2 \right).
\]
On the other hand, the orthogonal free Appell polynomials are only the Chebyshev polynomials of the second kind, with the ordinary generating function
\[
\frac{1}{1 - \left(x z - z^2 \right)}.
\]
Moreover, $R(z) = \alpha z^2$ appears naturally in combinatorial proofs of the usual, free, and other central-limit-type theorems (see for example Lecture~8 of \cite{Nica-Speicher-book}). Thus it is reasonable to consider $F\left(x z - \alpha z^2 \right)$ first, which we do in Theorem~\ref{Thm}. Its conclusion below indicates that interesting generating functions (and also, potentially, interesting non-commutative probability theories) arise precisely for $F$ covered by the Boas-Buck theorem above, plus in the exceptional case $F(z) = 1 + \log \frac{1}{1-z}$ not covered by that theorem. Next, we return to the question of generating functions of the more general form \eqref{Exp}, but only for $F$ of the special form just mentioned. Orthogonal polynomials with generating functions of this form are known, again with the exception of the special case $F(z) = 1 + \log \frac{1}{1-z}$. The conclusion is that such Meixner-type families are described by $4$ parameters for $\lambda = \infty$ (i.e. $F(z) = e^z$) and $\lambda = 1$, by $3$ parameters for other positive values of $\lambda$, and by $2$ parameters for $\lambda = 0$ (i.e. $F(z) = 1 + \log \frac{1}{1-z}$). See Remark~\ref{Remark:Meixner} and Proposition~\ref{Prop} below. So from the point of view of orthogonal polynomials, it is unlikely that there are non-commutative probability theories for other values of $\lambda$ strongly parallel to the classical and free theories.

After the article was submitted for publication, we learned from Professor Ben Cheikh that Question~\ref{Question:Main} has in fact been considered in the literature, in \cite{Al-Salam-Ultraspherical,von-Bachhaus,Ben-Cheikh-d-orthogonal-classical} (I am grateful to him for this information). Surprisingly, the characterization as described by Al-Salam included the ultraspherical and the limiting case of the Hermite polynomials, but seems to have missed the other limiting case of the Chebyshev polynomials of the first kind.

\textbf{Acknowledgments.} I am grateful to Harold Boas for many valuable comments, including Remark~\ref{Remark:Shift}. I would also like to thank the referee for the suggestion to use the Fa{\`a} di Bruno formula in the proof of Theorem~\ref{Thm}.

\section{The results}

\begin{Thm}
\label{Thm}
Let $\alpha > 0$ and $F(z) = \sum_{n=0}^\infty c_n z^n$ be a formal power series with $c_0 = 1$, $c_1 = c \neq 0$. Define the polynomials $\set{P_n : n \geq 0}$ via
\begin{equation}
\label{Eq:Generating function}
F \left(x z - \alpha z^2 \right) = \sum_{n=0}^\infty c_n P_n(x) z^n
\end{equation}
(if $c_n = 0$, $P_n$ is undefined). These polynomials form an orthogonal polynomial family (which is automatically monic) if and only if
\begin{itemize}
\item
$\set{P_n}$ are re-scaled ultraspherical polynomials,
\[
P_n(x) = C^{(\lambda)}_n \left( \sqrt{\frac{b}{4 \alpha}}\ x \right),
\]
for $\lambda > - \frac{1}{2}$, $\lambda \neq 0$, and $b > 0$. In this case
\[
F(z) = 1 + \frac{c}{\lambda b} \left( \frac{1}{(1 - b z)^{\lambda}} - 1 \right)
\]
for $c \neq 0$. The choice $c = \lambda$, $b=1$ gives simply $F(z) = \frac{1}{(1 - z)^{\lambda}}$.
\item
$\set{P_n}$ are re-scaled Chebyshev polynomials of the first kind,
\[
P_n(x) = T_n \left( \sqrt{\frac{b}{4 \alpha}}\ x \right),
\]
for $b > 0$. In this case
\[
F(z) = 1 + \frac{c}{b} \ln \left( \frac{1}{1 - b z} \right)
\]
for $c \neq 0$. The choice $c = b = 1$ gives simply $F(z) = 1 +  \ln \left( \frac{1}{1 - z} \right)$.
\item
$\set{P_n}$ are re-scaled Hermite polynomials,
\[
P_n(x) = H_n \left( \sqrt{\frac{a}{2 \alpha}}\ x \right)
\]
for $a > 0$. In this case
\[
F(z) = 1 + \frac{c}{a} \left( e^{a z} - 1 \right)
\]
for $c \neq 0$. The choice $c = a = 1$ gives simply $F(z) = e^{z}$.
\end{itemize}
\end{Thm}

\begin{Remark}
\label{Remark:Shift}
If equation~\eqref{Eq:Generating function} holds, so that
\begin{equation*}
1 + \sum_{n=1}^\infty c_n \left(x z - \alpha z^2 \right)^n = F \left(x z - \alpha z^2 \right) = 1 + \sum_{n=1}^\infty c_n P_n(x) z^n,
\end{equation*}
then clearly also
\begin{equation*}
1 + \sum_{n=1}^\infty C c_n \left(x z - \alpha z^2 \right)^n = F_C \left(x z - \alpha z^2 \right) = 1 + \sum_{n=1}^\infty C c_n P_n(x) z^n
\end{equation*}
for any $C \neq 0$ and $F_C(z) = 1 + C (F(z) - 1)$. This is the source of the free parameter $c$ in the theorem.
\end{Remark}

\begin{Remark}
\label{Remark:Meixner}
We now consider all orthogonal polynomials with more general generating functions of the form
\[
F \left(x w(z) - R(z) \right) = \sum_{n=0}^\infty c_n P_n(x) z^n,
\]
for $F$ as in the conclusion of the preceding theorem. For $F(z) = e^z$, corresponding to $\lambda = \infty$, this is a classical $4$-parameter family of Meixner polynomials. For $F(z) = \frac{1}{1-z}$, corresponding to $\lambda = 1$, these are the free Meixner polynomials, also a $4$-parameter family. For other positive $\lambda$, the orthogonal polynomials with generating functions of this form have been classified in \cite{Al-Salam-Verma,Demni-Ultraspherical}, and form only a $3$-parameter family. In the next proposition, we show that in the exceptional case $F(z) = 1 + \ln \frac{1}{1-z}$, corresponding to $\lambda = 0$, one only gets a $2$-parameter family. Note that we may re-write
\[
1 + \ln \frac{1}{1 - x w(z) + R(z)} = u(z) - \log \left( 1 - v(z) x \right).
\]
\end{Remark}

\begin{Prop}
\label{Prop}
The only orthogonal polynomials with the generating function of the form
\[
H(x,z) = u(z) - \log \left( 1 - v(z) x \right) = \sum_{n=0}^\infty \frac{1}{n} P_n(x) z^n,
\]
are Chebyshev polynomials of the first kind. Here
\[
u(z) = 1 + \sum_{n=1}^\infty \frac{1}{n} u_n z^n, \quad v(z) = \sum_{n=1}^\infty v_n z^n
\]
are formal power series.
\end{Prop}

\begin{proof}[Proof of the Theorem]
Using the Fa{\`a} di Bruno formula (see, for example, Theorem 3.3 in \cite{Aigner-Course-Enumeration}) to compute the coefficients of the composition of power series, or simply using the binomial expansion,
\[
\begin{split}
F \left(x z - \alpha z^2 \right)
& = \sum_{n=0}^\infty c_n \left(x z - \alpha z^2 \right)^n \\
& = \sum_{n=0}^\infty \sum_{0 \leq k \leq n/2} c_{n - k} \binom{n - k}{k} x^{n - 2 k} (- \alpha)^{k} z^n \\
\end{split}
\]
and so
\begin{equation}
\label{One}
c_n P_n(x) = \sum_{0 \leq k \leq n/2} c_{n - k} \binom{n - k}{k} (- \alpha)^{k} x^{n - 2 k}.
\end{equation}
If for some $n \geq 2$, $c_n = 0$ while $c_{n-1} \neq 0$, then comparing the coefficients of $x^{n-2}$ in equation~\eqref{One}, we see that the coefficient is zero on the left and non-zero on the right. So all $c_n \neq 0$, and we may denote $d_n = \frac{c_n}{c_{n-1}}$. Equation~\eqref{One} then becomes
\[
P_0(x) = 1, \quad P_1(x) = x, \quad P_2(x) = x^2 - \alpha d_2^{-1}, \quad P_3(x) = x^3 - 2 \alpha d_3^{-1} x
\]
and for $n \geq 4$,
\begin{equation}
\label{Two}
P_n(x) = x^n - \alpha (n-1) d_n^{-1} x^{n-2} + \alpha^2 \frac{(n-2)(n-3)}{2} d_n^{-1} d_{n-1}^{-1} x^{n-4} - \ldots.
\end{equation}
Using this equation for both $n$ and $n+1$, it follows that the leading term in $x P_n(x) - P_{n+1}(x)$ is
\[
\Bigl(\alpha n d_{n+1}^{-1} - \alpha (n-1) d_n^{-1} \Bigr) x^{n-1}.
\]
If we want the polynomials to be orthogonal, by Darboux-Favard-Stone theorem \cite{Chihara-book} they have to satisfy a three-term recursion relation
\[
x P_n = P_{n+1} + \beta_n P_n + \omega_n P_{n-1}
\]
(note that $\set{P_n}$ are clearly monic). We see that $\beta_n = 0$, and
\[
\omega_n = \alpha \Bigl(n d_{n+1}^{-1} - (n-1) d_n^{-1} \Bigr)
\]
for $n \geq 1$. Using \eqref{Two} again, it follows that for $n \geq 3$, the leading term in $x P_n(x) - P_{n+1}(x) - \omega_n P_{n-1}(x)$ is
\[
\biggl( \alpha^2 \frac{(n-2)(n-3)}{2} d_n^{-1} d_{n-1}^{-1} - \alpha^2 \frac{(n-1)(n-2)}{2} d_{n+1}^{-1} d_{n}^{-1}
- \alpha \Bigl((n-1) d_n^{-1} - n d_{n+1}^{-1}\Bigr) \alpha (n-2) d_{n-1}^{-1} \biggr) x^{n-3}.
\]
For this to be zero we need
\[
\frac{(n-3)}{2} d_n^{-1} d_{n-1}^{-1} - \frac{(n-1)}{2} d_{n+1}^{-1} d_{n}^{-1} - \Bigl((n-1) d_n^{-1} - n d_{n+1}^{-1} \Bigr) d_{n-1}^{-1} = 0,
\]
or
\[
\frac{(n-3)}{2} d_{n+1} - \frac{(n-1)}{2} d_{n-1} - \Bigl((n-1) d_{n+1} - n d_{n} \Bigr) = 0.
\]
Thus for $n \geq 3$,
\[
(n+1) d_{n+1} = 2 n d_n - (n-1) d_{n-1}.
\]
The general solution of this recursion is
\[
n d_n = a + b (n-1)
\]
for $n \geq 2$. Since all $d_n \neq 0$, $a, b$ cannot both be zero. Therefore
\[
c_n = \frac{a + b (n-1)}{n} c_{n-1} = \frac{\prod_{i=1}^{n-1} (a + i b)}{n!} c_1 = \frac{\prod_{i=1}^{n-1} (a + i b)}{n!} c
\]
for $n \geq 2$ and
\[
F(z) = 1 + c z + c \sum_{n=2}^\infty \frac{\prod_{i=1}^{n-1} (a + i b)}{n!} z^n.
\]
If $a \neq 0$, $b \neq 0$, then
\begin{equation}
\label{GF:ultraspherical}
\begin{split}
F(z)
& = 1 + c \sum_{n=1}^\infty \frac{\prod_{i=0}^{n-1} (-a/b - i)}{a n!} (-b z)^n \\
& = 1 + \frac{c}{a} \left( (1 - b z)^{-a/b} - 1 \right)
= 1 + \frac{c}{a} \left( \frac{1}{(1 - b z)^{a/b}} - 1 \right).
\end{split}
\end{equation}
If $a = 0$, $b \neq 0$, then
\begin{equation}
\label{GF:Chebyshev}
F(z) = 1 + c \sum_{n=1}^\infty \frac{b^{n-1}}{n} z^n
= 1 - \frac{c}{b} \ln(1 - b z)
= 1 + \frac{c}{b} \ln \left( \frac{1}{1 - b z} \right),
\end{equation}
which can also be obtained from the preceding formula by using L'H{\^o}pital's rule. Finally, if $a \neq 0$, $b=0$, then
\begin{equation}
\label{GF:Hermite}
F(z) = 1 + c \sum_{n=1}^\infty \frac{a^{n-1}}{n!} z^n = 1 + \frac{c}{a} \left( e^{a z} - 1 \right).
\end{equation}
Moreover,
\[
\omega_n = \alpha n \frac{(n-1) b + 2 a}{((n-1) b + a) (n b + a)}.
\]
Since for orthogonality, we need $\omega_n \geq 0$, clearly $b \geq 0$. If $b = 0$, then
\[
\omega_n = \frac{2 \alpha}{a} n > 0
\]
as long as $a > 0$. The polynomials with this recursion are the re-scaled Hermite polynomials. We recall \cite[Section 9.15]{Koekoek-S-book} that the generating function for standard (monic) Hermite polynomials is
\[
\sum_{n=0}^\infty \frac{1}{n!} H_n(x) z^n = \exp \left(x z - z^2/2 \right),
\]
which is of the form~\eqref{Eq:Generating function} with $F(z) = e^z$, and the generating function~\eqref{GF:Hermite} is obtained from it by a re-scaling and a shift from Remark~\ref{Remark:Shift}.

If $b > 0$, $a \neq 0$, we denote $\lambda = a/b$ and get
\[
\omega_n = \frac{\alpha}{b} \frac{n (n + 2 \lambda - 1)}{(n + \lambda - 1) (n + \lambda)}.
\]
Since
\[
\omega_1 = \frac{\alpha}{b} \frac{2}{1 + \lambda},
\]
we have $\lambda > -1$. Since
\[
\omega_2 = \frac{\alpha}{b} \frac{2 (1 + 2 \lambda)}{(1 + \lambda) (2 + \lambda)},
\]
we have moreover $\lambda > - \frac{1}{2}$. It is now easy to see that this condition suffices for the positivity of all $\omega_n$; indeed, the corresponding polynomials are the re-scaled ultraspherical polynomials. We recall \cite[Section 9.8.1]{Koekoek-S-book} that the generating function for standard (monic) ultraspherical  polynomials is
\[
\sum_{n=0}^\infty \frac{2^n \prod_{i=0}^{n-1} (\lambda - i)}{n!} C^{(\lambda)}_n(x) z^n = \frac{1}{\left(1 - 2 x z + z^2 \right)^{\lambda}},
\]
which is of the form~\eqref{Eq:Generating function} with $F(z) = \frac{1}{(1 - 2 z)^{\lambda}}$, and the generating function~\eqref{GF:ultraspherical} is obtained from it by a re-scaling and a shift from Remark~\ref{Remark:Shift}.

Finally, if $b > 0$, $\lambda = a = 0$, then
\[
\omega_n = \frac{\alpha}{b}
\]
for $n \geq 2$, but $\omega_1 = 2 \frac{\alpha}{b}$. These are precisely recursion coefficients for the re-scaled Chebyshev polynomials of the fist kind. The standard generating function \cite[Section 9.8.2]{Koekoek-S-book} for (monic) Chebyshev polynomials of the first kind is
\begin{equation*}
1 + \sum_{n=1}^\infty 2^{n-1} T_n(x) z^n = \frac{1 - x z}{1 - 2 x z + z^2},
\end{equation*}
so it is not of the form~\eqref{Eq:Generating function}. However,
\begin{equation*}
\sum_{n=1}^\infty 2^{n-1} T_n(x) z^{n-1} = \frac{1}{z} \left( \frac{1 - x z}{1 - 2 x z + z^2} - 1 \right)
= \frac{x - z}{1 - 2 x z + z^2}.
\end{equation*}
Term-by-term integration with respect to $z$ gives
\[
C + \sum_{n=1}^\infty \frac{2^{n-1}}{n} T_n(x) z^{n}
= 1 - \frac{1}{2} \ln \left(1 - 2 x z + z^2 \right)
= 1 + \frac{1}{2} \ln \left( \frac{1}{1 - 2 \left(x z - z^2/2 \right)} \right).
\]
with $C = 1$, which is of the form~\eqref{Eq:Generating function} with $F(z) = 1 + \frac{1}{2} \ln \left( \frac{1}{1 - 2 z} \right)$. The generating function~\eqref{GF:Chebyshev} is obtained from it by a re-scaling and a shift from Remark~\ref{Remark:Shift}.
\end{proof}

\begin{proof}[Proof of the Proposition]
We use the same method as in the theorem. By assumption,
\[
1 + \sum_{n=1}^\infty \frac{1}{n} u_n z^n + \sum_{k=1}^\infty \frac{1}{k} v(z)^k x^k = \sum_{n=0}^\infty \frac{1}{n} P_n(x) z^n.
\]
The left-hand-side is
\[
1 + (u_1 + v_1 x) z + (\frac{1}{2} u_2 + v_2 x + \frac{1}{2} v_1^2 x^2) z^2 + \ldots.
\]
By applying an affine transformation to the measure of orthogonality, we may assume that it has mean zero and unit variance. Thus $P_0(x) = 1$,
\[
P_1(x) = v_1 x + u_1 = x,
\]
so that $v_1 = 1$, $u_1 = 0$, and
\[
P_2(x) = x^2 + 2 v_2 x + 2 u_2,
\]
so $u_2 = - 1$. Expanding further, the left-hand-side is
\[
1 + P_1(x) z + \frac{1}{2} P_2(x) z^2 + \sum_{n=3}^\infty \left( \frac{1}{n} x^n + v_2 x^{n-1} + (v_3 + \frac{n-3}{2} v_2^2) x^{n-2} + \ldots + \frac{1}{n} u_n \right),
\]
so for $n \geq 2$,
\[
P_n(x) = x^n + n v_2 x^{n-1} + n (v_3 + \frac{n-3}{2} v_2^2) x^{n-2} + \ldots + u_n.
\]
On the other hand,
\[
x P_n = P_{n+1} + \beta_n P_n + \omega_n P_{n-1}.
\]
Thus
\begin{gather*}
n v_2 = (n+1) v_2 + \beta_n, \\
n (v_3 + \frac{n-3}{2} v_2^2) = (n + 1) (v_3 + \frac{n-2}{2} v_2^2) + \beta_n n v_2 + \omega_n
\end{gather*}
and
\[
0 = u_{n+1} + \beta_n u_n + \omega_n u_{n-1}.
\]
It follows that
\[
\beta_n = - v_2, \quad \beta_0 = 0,
\]
\[
\omega_n = v_2^2 - v_3, \quad \omega_1 = 1,
\]
\[
0 = u_{n+1} - v_2 u_n + (v_2^2 - v_3) u_{n-1}, \quad u_1 = 0, u_2 = -1.
\]
Thus
\[
0 = z + u'(z) - v_2 z u'(z) + (v_2^2 - v_3) z^2 u'(z),
\]
and
\begin{equation}
\label{u-prime}
u'(z) = - \frac{z}{1 - v_2 z + (v_2^2 - v_3) z^2}.
\end{equation}
It follows from the recursion formulas that the polynomials belong to the normalized free Meixner class. Their ordinary generating function is known \cite{Al-Salam-Verma,AnsMeixner}, and has the form
\[
\sum_{n=0}^\infty P_n(x) z^n = \frac{p(z)}{q(z) - z x},
\]
where $p$ and $q$ are quadratic polynomials with constant term $1$. Thus
\[
z \partial_z H(x,z) + 1 = 1 + z u'(z) + \frac{z v'(z) x}{1 - v(z) x} = \frac{p(z)}{q(z) - z x} = \sum_{n=0}^\infty P_n(x) z^n,
\]
from which we deduce after a short computation that $v(z) = z e^{u(z)}$, $q(z) = e^{-u(z)}$, $p(z) = e^{-u(z)} (1 + z u'(z))$. In particular,
\[
u'(z) = - \frac{q'(z)}{q(z)}.
\]
Comparing with equation~\eqref{u-prime}, we conclude that $v_2 = 0$ and $v_3 = - \frac{1}{2}$. Thus $\beta_n = 0$, $\omega_1 = 1$, and $\omega_n = \frac{1}{2}$ for $n > 1$. These are the Jacobi parameters for the Chebyshev polynomials of the first kind. As observed in the preceding theorem, they have the generating function of the desired form.
\end{proof}


\def\cprime{$'$}
\providecommand{\bysame}{\leavevmode\hbox to3em{\hrulefill}\thinspace}
\providecommand{\MR}{\relax\ifhmode\unskip\space\fi MR }
\providecommand{\MRhref}[2]{%
  \href{http://www.ams.org/mathscinet-getitem?mr=#1}{#2}
}
\providecommand{\href}[2]{#2}

\end{document}